# THE STRUCTURE OF DGA RESOLUTIONS OF MONOMIAL IDEALS

LUKAS KATTHÄN

ABSTRACT. Let $I \subseteq \Bbbk[x_1, \dots, x_n]$ be a squarefree monomial ideal in a polynomial ring. In this paper we study multiplications on the minimal free resolution $\mathbb{F}$ of $\Bbbk[x_1, \dots, x_n]/I$. In particular, we characterize the possible vectors of total Betti numbers for such ideals which admit a differential graded algebra (DGA) structure on $\mathbb{F}$. We also show that under these assumptions the maximal shifts of the graded Betti numbers are subadditive.

On the other hand, we present an example of a strongly generic monomial ideal which does not admit a DGA structure on its minimal free resolution. In particular, this demonstrates that the Hull resolution and the Lyubeznik resolution do not admit DGA structures in general.

Finally, we show that it is enough to modify the last map of $\mathbb{F}$ to ensure that it admits the structure of a DG algebra.

## INTRODUCTION

Let $S$ be a polynomial ring over a field, $I \subseteq S$ be a monomial ideal, and let $\mathbb{F}$ denote the minimal free resolution of $S/I$ over $S$. The multiplication map $S/I \otimes_S S/I \to S/I$ can be lifted to a "multiplication" $* : \mathbb{F} \otimes_S \mathbb{F} \to \mathbb{F}$, which in general is associative and graded-commutative only up to homotopy. Moreover, $*$ is unique only up to homotopy. It is known that $*$ can always chosen to be graded-commutative "on the nose" (cf. [BE77, Prop 1.1]), but in general it is not possible to choose $*$ such that it is associative (cf. [Avr74, Appendix], [Avr81]). If $*$ is graded-commutative and associative, then it gives $\mathbb{F}$ the structure of a differential graded algebra (DGA). In this situation, we say that $S/I$ admits a *minimal DGA resolution*. Starting with the work of Buchsbaum and Eisenbud [BE77], a lot of research has been devoted to the study of DGA resolutions, see for example [Avr81] or [Kus94] and the references therein. In the present paper, we study the case where $I$ is generated by squarefree monomials and consider only multiplications which respect the multigrading. These restrictions allow us to prove much stronger statements than in the general situation. Our results can be organized along the following questions:

(1) Which ideals admit a minimal DGA resolution? Which constructions of (not necessarily minimal) resolutions yield DGA resolutions?
(2) What are the consequences if a given ideal admits a minimal DGA resolution?
(3) How unique or non-unique are the DGA structures?
(4) What can be said about the structure of DGA resolutions as algebras?

There are a few classes of ideals in local rings which are known to admit minimal DGA resolutions, see [Kus94] for an overview. In the context of monomial ideals, this is the case for stable ideals [Pee96], matroidal ideals [Skö11] and edge ideals of cointerval graphs [Skö16]. Another class of ideals whose minimal free resolutions are well-understood are generic monomial ideals. However, in Theorem 5.1, we present an example of a (strongly) generic monomial ideal which does not admit a minimal DGA resolution. The example also demonstrates that the hull resolution [BS98] and the Lyubeznik resolution [Lyu88] do not admit a DGA structure in general.

On the positive side, we show in Theorem 2.1 that for every monomial ideal $I$ there always exists a monomial $m \in S$ such that $S/(mI)$ admits a minimal DGA resolution. Note that the minimal free resolutions of $S/I$ and $S/(mI)$ differ only in the last map. Therefore, one can interpret this result as saying that the other maps in the resolution do not contain obstructions to the existence of a DGA structure. This result actually holds in greater generality for homogeneous ideals in $S$ and even for ideals in regular local domains.

Concerning the second question, it is a beautiful result by Buchsbaum and Eisenbud [BE77, Proposition 1.4], that if $I$ admits a minimal DGA resolution, then one can conclude certain lower bounds on the Betti numbers of $S/I$. In the context of squarefree monomial ideals, we sharpen this result by completely characterizing the possible total Betti numbers of ideals admitting a minimal DGA resolution (Theorem 4.1). Moreover, we show that if $I$ is such an ideal, then its graded Betti numbers satisfy the subadditivity condition of [ACI15] (Proposition 4.4).

It is well-known that DGA structures on minimal free resolutions do not need to be unique. Our Example 3.2 demonstrates that this even fails for squarefree monomial ideals, thus answering a question of Buchsbaum and Eisenbud [Pee10, Open Problem 31.4]. On the other hand, we show that at least a certain part of the multiplication is indeed unique, see Proposition 3.1 for the precise statement.

Concerning the fourth question, we show that $\mathbb{F}$ is essentially generated in degree 1 with respect to any multiplication, cf. Proposition 3.4. As a consequence, whenever $I$ admits a minimal DGA resolution $\mathbb{F}$, then it is a quotient of the Taylor resolution (which has a canonical DGA structure) by a DG-ideal, see Theorem 3.6. Thus one can study the possible DGA structures on $\mathbb{F}$ by considering DG-ideals in the Taylor resolution.

The key idea behind most of our results is the following simple observation: $\mathbb{F}$ is a free $S$-module which is generated in squarefree degrees, so for any homogeneous element $f \in \mathbb{F}$, we can find an element $f' \in \mathbb{F}$ of squarefree degree, such that $f = mf'$ for a monomial $m \in S$. We call $f'$ the *squarefree part* of $f$.

This paper is structured as follows. In Section 1 we recall some preliminaries about multiplicative structures on resolution, the Taylor resolution and the Scarf complex, and we introduce the squarefree part of an element $f \in \mathbb{F}$. In the next section we prove Theorem 2.1 about the existence of minimal DGA resolutions after multiplication with an element. The following Section 3 is devoted to the study of the structure of multiplication on $\mathbb{F}$. After that, in Section 4 we derive two consequences of the existence of a minimal DGA resolution, namely our characterization of the total Betti numbers and the subadditivity of syzygies.

In the penultimate Section 5 we present an example of a generic monomial ideal which does not admit a minimal DGA resolution, and in the last section we study multiplication on resolutions of ideals which are not necessarily squarefree.

## 1. Preliminaries and notation

Throughout the paper, let $\Bbbk$ be a field and $S = \Bbbk[x_1, \ldots, x_n]$ be a polynomial ring over $\Bbbk$, endowed with the fine $\mathbb{Z}^n$-grading. We denote by $\mathfrak{m} := \langle x_1, \ldots, x_n \rangle$ the maximal homogeneous ideal of $S$. Further, let $I \subseteq S$ be a monomial ideal and let
$$\mathbb{F}: \quad \cdots \xrightarrow{\partial} \mathbb{F}_2 \xrightarrow{\partial} \mathbb{F}_1 \xrightarrow{\partial} \mathbb{F}_0 \to S/I \to 0$$
denote the minimal free resolution of $S/I$. Note that $\mathbb{F}$ is $\mathbb{Z}^n$-graded and the differential $\partial$ respects the multigrading. The multigraded Betti numbers of $S/I$ are denoted by $\beta^S_{i,\mathbf{a}}(S/I) = \dim_{\Bbbk} \operatorname{Tor}^S_i(S/I, \Bbbk)_\mathbf{a}$. We often identify $\mathbb{F}$ with the free multigraded $S$-module $\bigoplus_i \mathbb{F}_i$. For an element $f \in \mathbb{F}$ we denote its multidegree by $\deg(f)$ and its homological degree by $|f|$.

We consider the componentwise order on $\mathbb{Z}^n$ and write $\mathbf{a} \wedge \mathbf{b}$ and $\mathbf{a} \vee \mathbf{b}$ for the componentwise minimum and maximum of $\mathbf{a}, \mathbf{b} \in \mathbb{Z}^n$, respectively.

### 1.1. Multiplications on resolutions.
We recall the definition of a DGA structure.

**Definition 1.1.** A *differential graded algebra* (DGA) structure on $\mathbb{F}$ is an $S$-linear map $* : \mathbb{F} \otimes_S \mathbb{F} \to \mathbb{F}$ satisfying the following axioms for $a, b, c \in \mathbb{F}$:
(1) $*$ extends the usual multiplication on $\mathbb{F}_0 = S$,
(2) $\partial(a * b) = \partial(a) * b + (-1)^{|a|} a * \partial b$ (Leibniz rule),
(3) $|a * b| = |a| + |b|$ (homogeneity with respect to the homological grading),
(4) $a * b = (-1)^{|a| \cdot |b|} b * a$ (graded commutativity), and
(5) $(a * b) * c = a * (b * c)$ (associativity).

We will also consider non-associative multiplications. To make this precise, we call a map $* : \mathbb{F} \otimes_S \mathbb{F} \to \mathbb{F}$ a *multiplication* if it satisfies all the axioms of a DGA except possibly the associativity. Moreover, we make the convention that in this paper, *every multiplication respects the multigrading on $\mathbb{F}$*, unless specified otherwise. The only occasion when we consider more general multiplication is in Section 5.

While $\mathbb{F}$ does not always admit a DGA structure, it always admits a multiplication, cf. [BE77, Prop 1.1], and the multiplication is unique up to homotopy. Explicitly, this means that when $*_1, *_2$ are two multiplications, then there exists a map $\sigma : \mathbb{F} \otimes_S \mathbb{F} \to \mathbb{F}$ raising the homological degree by 1 such that $a *_1 b = a *_2 b + \partial \sigma(a, b) + \sigma(\partial a, b) + (-1)^{|a|} \sigma(a, \partial b)$ for $a, b \in \mathbb{F}$.

### 1.2. The Taylor resolution and the Scarf complex.
We recall the definitions of the Taylor resolution and the Scarf complex of $I$. We refer the reader to p.67 and Section 6.2 of [MS05], respectively, for further information about these constructions.

Let $G(I)$ denote the set of minimal monomial generators of $I$ and choose a total order $\prec$ on $G(I)$. The choice of the order affects only the signs in the computations.

**Definition 1.2.** The *Taylor resolution* $\mathbb{T}$ of $S/I$ is the complex of free $S$-modules with basis $\{\mathbf{g}_W : W \subseteq G(I)\}$. The basis elements are graded by $|\mathbf{g}_W| := \#W$ and $\deg \mathbf{g}_W := \deg m_W$, where $m_W := \operatorname{lcm}(m : m \in W)$ for $W \subseteq G(I)$. Further, the differential is given by

$$\partial \mathbf{g}_W = \sum_{m \in W} (-1)^{\sigma(m,W)} \frac{m_W}{m_{W \setminus \{m\}}} \mathbf{g}_{W \setminus \{m\}},$$

where $\sigma(m, W) := \#\{m' \in W, m' \prec m\}$.

The Taylor resolution is a free resolution of $S/I$, but typically not a minimal one. It was shown by Gemeda [Gem76] (see also [Pee10, Proposition 31.3]) that it carries a DGA structure with the multiplication given by

$$\mathbf{g}_W * \mathbf{g}_V = \begin{cases} (-1)^{\sigma(W,V)} \frac{m_W m_V}{m_{W \cup V}} \mathbf{g}_{W \cup V} & \text{if } W \cap V = \emptyset, \\ 0 & \text{otherwise,} \end{cases}$$

where $\sigma(W, V) := \#\{(m, m') \in W \times V : m' \prec m\}$. To simplify the notation we occasionally write $\mathbf{g}_{m_1 m_2 \dots}$ instead of $\mathbf{g}_{\{m_1, m_2, \dots\}}$ for $m_1, m_2, \dots \in G(I)$.

**Definition 1.3.** The *Scarf complex* $\Delta_I$ of $S/I$ is the simplicial complex with vertex set $G(I)$, where a set $W \subseteq G(I)$ is contained in $\Delta_I$ if and only if there does not exist another set $V \subseteq G(I)$ with $V \neq W$ and $\operatorname{lcm}(m : m \in W) = \operatorname{lcm}(m : m \in V)$.

Moreover, the *algebraic Scarf complex* $\mathcal{F}_{\Delta_I}$ is the subcomplex of the Taylor resolution generated by the generators $\mathbf{g}_W, W \in \Delta_I$.

The algebraic Scarf complex is always a subcomplex of the minimal free resolution of $S/I$ (cf. [MS05, Proposition 6.12]), but in general it is not acyclic. We will therefore use the notation $\mathbf{g}_W$ for generators of $\mathbb{F}$ which lie in the Scarf complex.

1.3. **The squarefree part.** The following very basic observation turns out to be the key to many of our results:

**Lemma 1.4.** *Let $F$ be a free $S$-module and $\mathbf{a} \in \mathbb{N}^n$. Assume that the degree of every generator of $F$ is less or equal to $\mathbf{a}$. Then every element $f \in F$ with $\deg f \not\leq \mathbf{a}$ can be written as $f = mf'$ for a monomial $m \in S$ and $f' \in F$ with $\deg f' = (\deg f) \wedge \mathbf{a}$. Both $m$ and $f'$ are uniquely determined.*

This lemma is almost obvious, but we include a proof for completeness.

*Proof.* Choose a basis of $F$ and let $G$ be the set of those basis elements whose degrees are less or equal than $\deg f$. By assumption, $f$ can be written as a linear combination of elements of $G$, and we write $\lambda_\mathbf{g}$ for the coefficient of $\mathbf{g}$ in $\Bbbk$ in this expansion. Then it holds that

$$f = \sum_{\mathbf{g} \in G} \lambda_\mathbf{g} x^{\deg(f) - \deg(\mathbf{g})} \mathbf{g} = x^{\deg(f) - (\mathbf{a} \wedge \deg(f))} \sum_{\mathbf{g} \in G} \lambda_\mathbf{g} x^{(\mathbf{a} \wedge \deg(f)) - \deg(\mathbf{g})} \mathbf{g}.$$

All exponents are non-negative by our hypothesis, so we have shown the existence of the claimed factorization.

The uniqueness of $m$ is trivial, because there is only one monomial of the correct multidegree. Using this, the uniqueness of $f'$ follows from the fact that $S$ is a domain and $F$ is torsion-free. □

**Definition 1.5.** *Let $F$ be a free $S$-module, such that the degree every generator is squarefree, i.e. less or equal than $(1, \ldots, 1) \in \mathbb{N}^n$. In the situation of Lemma 1.4, we call $f'$ the* squarefree part *of $f$ and denote it by $|f|_{\mathrm{sqf}}$.*

Note that the map $f \mapsto |f|_{\mathrm{sqf}}$ is $\Bbbk$-linear, but not $S$-linear. We will apply this definition almost exclusively to minimal free resolutions $\mathbb{F}$ of quotients $S/I$ by squarefree monomial ideals $I \subseteq S$.

## 2. Existence of DGA structures up to multiplication with an element

Our first result is the following theorem. It implies that it is enough to modify the last map of a minimal free resolution to ensure the existence of a DGA structure.

**Theorem 2.1.** *Let $S$ be a regular local ring and $I \subseteq S$ an ideal. There exists an element $s \in S, s \neq 0$ such that the minimal free resolution of $S/(sI)$ admits a DGA structure.*

*The same conclusion holds if $S$ is a polynomial ring and $I$ is a homogeneous ideal. In this case, $s$ can be chosen homogeneous. If $I$ is even a monomial ideal (and $S$ a polynomial ring), then $s$ can be chosen to be the least common multiple of the generators of $I$.*

We need the following lemma for the proof of the Theorem.

**Lemma 2.2.** *Let $Q = \Bbbk[t_1^{\pm}, \ldots, t_n^{\pm}]$ be a $\mathbb{Z}^n$-graded Laurent polynomial ring for $n \geq 0$. Let*
$$\mathbb{F}: \quad 0 \to \mathbb{F}_p \to \cdots \to \mathbb{F}_1 \to \mathbb{F}_0 \to 0$$
*be an exact complex of graded $Q$-modules and assume that $\mathbb{F}_0 = Q$. Then the multiplication on $\mathbb{F}_0$ can be extended to a graded-commutative DGA structure on $\mathbb{F}$.*

*Proof.* Let $\partial$ denote the differential of $\mathbb{F}$. We claim that there exists a map $\sigma : \mathbb{F} \to \mathbb{F}$ of homological degree 1 such that $\partial \circ \sigma + \sigma \circ \partial = \mathrm{id}_{\mathbb{F}}$ and $\sigma \circ \sigma = 0$. Indeed, every graded module over $Q$ is free (cf. [GW78, Theorem 1.1.4]), so we can choose splittings $\mathbb{F}_i \cong V_i \oplus \partial(\mathbb{F}_{i+1})$. Note that $\partial(\mathbb{F}_{i+1}) = \partial(V_{i+1})$ and the restriction of $\partial$ to $V_{i+1}$ is injective. We define $\sigma_i|_{V_i} = 0$ and $\sigma_i(\partial(f)) = f$ for $f \in V_{i+1}$. It is not difficult to see that this gives indeed a map $\sigma$ with the claimed properties.

We define the DGA structure on $\mathbb{F}$ inductively by the formula
$$a * b := \begin{cases} ab & \text{if } |a| = 0 \text{ or } |b| = 0 \\ \sigma\left(\partial(a) * b + (-1)^{|a|} a * \partial(b)\right) & \text{otherwise,} \end{cases}$$
where the multiplication in the first case is the one from $\mathbb{F}_0 = Q$.

This multiplication clearly satisfies the Leibniz rule and it extends the multiplication on $\mathbb{F}_0$. It remains to show that $*$ is graded-commutative and associative. Note that $\sigma = \sigma \circ \partial \circ \sigma$. We proceed by induction on the homological degree, the base case being clear. It holds for $a, b, c \in \mathbb{F}$:

$$\begin{aligned}
a * (b * c) &= \sigma\big(\partial(a * (b * c))\big) \\
&= \sigma\big(\partial a * (b * c) + (-1)^{|a|} a * (\partial b * c) + (-1)^{|a|+|b|} a * (b * \partial c)\big) \\
&\stackrel{(\#)}{=} \sigma\big((\partial a * b) * c + (-1)^{|a|}(a * \partial b) * c + (-1)^{|a|+|b|}(a * b) * \partial c\big) \\
&= \sigma\big(\partial((a * b) * c)\big) \\
&= (a * b) * c
\end{aligned}$$

where in $(\#)$ we use the induction hypothesis. The commutativity is verified analogously. $\square$

*Proof of Theorem 2.1.* Let $\mathbb{F}$ be the minimal free resolution of $S/I$. In the local case, let $Q$ be the field of fractions of $I$, while in the polynomial case, let $Q$ be the subring of the field of fractions of $S$ where we adjoin inverses for all *homogeneous* elements of $S$. Set $\mathbb{F}_Q := \mathbb{F} \otimes_S Q$. In both cases, $Q$ is a multivariate Laurent polynomial ring over some field $\Bbbk$, though without variables in the local case. Moreover, $Q$ is flat over $S$ and hence $\mathbb{F}_Q$ is exact. So by Lemma 2.2, it can be endowed with a DGA structure $*$.

Choose a basis for each $\mathbb{F}_i$. Then the multiplication on $\mathbb{F}_Q$ can be represented as a matrix with entries in $Q$. Hence we can choose an element $s \in S$ such that $sq \in S$ for each entry $q$ of this matrix. Let $\mathbb{F}'$ be the subcomplex of $\mathbb{F}$ defined by $\mathbb{F}'_0 := \mathbb{F}_0$ and $\mathbb{F}'_i := s\mathbb{F}_i$ for $i \geq 1$. We claim that $\mathbb{F}'$ is closed under multiplication. Indeed, for $sa, sb \in \mathbb{F}'$ the choice of $s$ implies that $s(a * b) \in \mathbb{F}$ and thus $(sa * sb) \in \mathbb{F}'$.

Note that $\mathbb{F}'$ is isomorphic to $\mathbb{F}$ except in degree 0, so in particular it is exact in every other degree. In degree 0, it holds that $H_0(\mathbb{F}') = \mathbb{F}_0/\partial(s\mathbb{F}_1) = S/(sI)$. Thus $\mathbb{F}'$ is the minimal free resolution of $S/(sI)$.

Finally, consider the case that $I$ is a monomial ideal in a polynomial ring $S = \Bbbk[x_1, \ldots, x_n]$. Let $a, b \in \mathbb{F}$ be two homogeneous elements. The product $a * b \in \mathbb{F}_Q$ can be written as a sum

$$\sum_{\mathbf{g} \in \mathcal{B}} \lambda_{\mathbf{g}} x^{\deg(a) + \deg(b) - \deg(\mathbf{g})} \mathbf{g}$$

where $\mathcal{B}$ is an $S$-basis of $\mathbb{F}$ and $\lambda_{\mathbf{g}} \in \Bbbk$. Let now $s$ be the lcm of all generators of $I$. Then the multidegrees of all elements of $\mathcal{B}$ are less or equal to $\deg(s)$, and hence $s(a * b) \in \mathbb{F}$. Now one can argue as above. $\square$

## 3. Properties of multiplications

3.1. **Uniqueness.** The multiplication on $\mathbb{F}$ is in general not unique. However, the part of the multiplication which lives on the algebraic Scarf complex is uniquely determined, as the next proposition shows.

**Proposition 3.1.** *Let $I \subseteq S$ be a squarefree monomial ideal with minimal free resolution $\mathbb{F}$, and let further $*$ be a multiplication on $\mathbb{F}$. Then, for $V, W \in \Delta_I$ with $V \cup W \in \Delta_I$, it holds that*

$$\mathbf{g}_W * \mathbf{g}_V = \begin{cases} (-1)^{\sigma(W,V)} \frac{m_W m_V}{m_{W \cup V}} \mathbf{g}_{W \cup V} & \text{if } W \cap V = \emptyset, \\ 0 & \text{otherwise}; \end{cases}$$

*where $\sigma(W, V)$ is defined as above. In particular, the product $\mathbf{g}_W * \mathbf{g}_V$ does not depend on the choice of $*$.*

*Proof.* Let $\mathcal{B}$ be an $S$-basis for $\mathbb{F}$. For every $\mathbf{g} \in \mathcal{B}$ there exists a (not necessarily unique) subset $U_\mathbf{g} \subseteq G(I)$ such that $\deg \mathbf{g} = \deg m_{U_\mathbf{g}}$ and $|\mathbf{g}| = \#U_\mathbf{g}$. This follows from the fact that $\mathbb{F}$ is a direct summand of the Taylor resolution.

Expand $\mathbf{g}_W * \mathbf{g}_V$ in the basis $\mathcal{B}$:
$$\mathbf{g}_W * \mathbf{g}_V = \sum_{\mathbf{g} \in \mathcal{B}} c_\mathbf{g} \mathbf{g}.$$

Assume that $c_\mathbf{g} \neq 0$ for some $\mathbf{g} \in \mathcal{B}$. As $I$ is squarefree, it holds that $\deg m_{U_\mathbf{g}} = \deg \mathbf{g} \leq \deg |\mathbf{g}_W * \mathbf{g}_V|_{\text{sqf}} = \deg(\mathbf{g}_W) \vee \deg(\mathbf{g}_V) = \deg m_{W \cup V}$. Hence $\deg m_{U_\mathbf{g} \cup W \cup V} = \deg m_{W \cup V}$. But $W \cup V$ is contained in the Scarf complex, and hence it follows that $U_\mathbf{g} \cup W \cup V = W \cup V$ and hence $U_\mathbf{g} \subseteq V \cup W$. On the other hand, it holds that $\#U_\mathbf{g} = |\mathbf{g}| = \#W + \#V \geq \#(W \cup V)$. So we can conclude that $U_\mathbf{g} = V \cup W$ and that $\#W + \#V = \#(W \cup V)$. Hence if $W$ and $V$ intersect non-trivially then $\mathbf{g}_W * \mathbf{g}_V$ is zero, and if they are disjoint, then $\mathbf{g}_W * \mathbf{g}_V = \lambda \frac{m_W m_V}{m_{W \cup V}} \mathbf{g}_{W \cup V}$ for a scalar $\lambda \in \Bbbk$. It is clear by induction that the scalars $\lambda$ are uniquely determined by the Leibniz rule. Moreover, the values given in the claim satisfy the Leibniz rule, hence we have determined the product $\mathbf{g}_W * \mathbf{g}_V$. □

The following example shows that the multiplication is not unique in general. This answers a question of Buchsbaum and Eisenbud [Pee10, Open Problems 31.4].

**Example 3.2.** Let $I \subseteq S = \mathbb{Q}[x_1, x_2, x_3, x_4, y, z]$ be the ideal with generators $a := x_1 x_2 y, b := x_2 x_3, c := x_3 x_4 z$ and $d := x_4 x_1$. One can verify with `Macaulay2` that the minimal free resolution $\mathbb{F}$ of $S/I$ equals its algebraic Scarf complex. The Scarf complex $\Delta_I$ is the following simplicial complex:

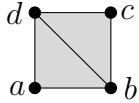

Consider the product $\mathbf{g}_a * \mathbf{g}_c$. Both $a$ and $c$ are vertices of $\Delta_I$, but $\{a, c\} \notin \Delta_I$, so the preceding Proposition does not apply. And indeed, for any $\lambda \in \Bbbk$ the choice
$$\mathbf{g}_a * \mathbf{g}_c := \lambda(x_4 z \mathbf{g}_{ab} + x_1 y \mathbf{g}_{bc}) + (1 - \lambda)(x_3 z \mathbf{g}_{ad} - x_2 y \mathbf{g}_{cd})$$
satisfies the Leibniz rule. This can be extended to a multiplication on $\mathbb{F}$, which is even a DGA structure, because the length of the resolution is three (cf. [BE77, Prop 1.3]).

The next example illustrates that one cannot omit the assumption of $I$ being squarefree from Proposition 3.1.

**Example 3.3.** Consider the ideal $I \subseteq S = \Bbbk[x, y, z]$ generated by $a := x^2, b := xy$ and $c := xz$. The algebraic Scarf complex of this ideal coincides with its Taylor resolution $\mathbb{T}$. So if $I$ were squarefree, then the DGA structure on $\mathbb{T}$ described above were the only possible multiplication. However, we can modify the multiplication by setting $\mathbf{g}_b * \mathbf{g}_c := z\mathbf{g}_{ab} - y\mathbf{g}_{ac}$. The Leibniz rule requires that in addition we set $\mathbf{g}_b * \mathbf{g}_{ac} := \mathbf{g}_c * \mathbf{g}_{ab} := 0$. See Table 1 for the full multiplication table of this new multiplication. It is even a DGA structure, again because the projective dimension of $S/I$ is less than four.

| $*$ | $\mathbf{g}_a$ | $\mathbf{g}_b$ | $\mathbf{g}_c$ | $\mathbf{g}_{ab}$ | $\mathbf{g}_{ac}$ | $\mathbf{g}_{bc}$ | $\mathbf{g}_{abc}$ |
|---|---|---|---|---|---|---|---|
| $\mathbf{g}_a$ | 0 | $x\mathbf{g}_{ab}$ | $x\mathbf{g}_{ac}$ | 0 | 0 | $x\mathbf{g}_{abc}$ | 0 |
| $\mathbf{g}_b$ | $-x\mathbf{g}_{ab}$ | 0 | $z\mathbf{g}_{ab} - y\mathbf{g}_{ac}$ | 0 | 0 | 0 | 0 |
| $\mathbf{g}_c$ | $-x\mathbf{g}_{ac}$ | $-z\mathbf{g}_{ab} + y\mathbf{g}_{ac}$ | 0 | 0 | 0 | 0 | 0 |
| $\mathbf{g}_{ab}$ | 0 | 0 | 0 | 0 | 0 | 0 | 0 |
| $\mathbf{g}_{ac}$ | 0 | 0 | 0 | 0 | 0 | 0 | 0 |
| $\mathbf{g}_{bc}$ | $x\mathbf{g}_{abc}$ | 0 | 0 | 0 | 0 | 0 | 0 |
| $\mathbf{g}_{abc}$ | 0 | 0 | 0 | 0 | 0 | 0 | 0 |

TABLE 1. The multiplication table of Example 3.3.

3.2. **Structure.** The next result shows that if $I$ is squarefree, then $\mathbb{F}$ is essentially generated in degree 1 for *any* multiplication on it:

**Proposition 3.4.** *Let $I \subseteq S$ be a squarefree monomial ideal with minimal free resolution $\mathbb{F}$. Let $*$ be a multiplication on $\mathbb{F}$. Then every homogeneous element $f \in \mathbb{F}$ of squarefree degree with $|f| \geq 1$ can be written as*
$$f = \sum_j |g_{1,j} * (g_{2,j} * (\cdots (g_{|f|-1,j} * g_{|f|,j}) \cdots ))|_{\mathrm{sqf}},$$
*for $g_{\ell,j} \in \mathbb{F}_1$ for $1 \leq \ell \leq |f|$ and all $j$.*

*Proof.* It is sufficient to show the claim for $f \in \mathbb{F} \setminus \mathfrak{m}\mathbb{F}$, because otherwise $f$ is a linear combination of such elements.

Let $g_{1,1} \in \mathbb{F}_1$ be a basis element with $\deg g_{1,1} \leq f$. Such an element exists, because $f$ is a linear combination of basis elements, and the degree of each basis element is an lcm of degrees of basis elements of $\mathbb{F}_1$. Note that $\deg g_{1,1} \leq \deg f$ implies that $f = |\partial(g_{1,1}) * f|_{\mathrm{sqf}}$. Hence it holds that
$$r := f - |g_{1,1} * \partial f|_{\mathrm{sqf}} = |\partial(g_{1,1} * f)|_{\mathrm{sqf}} = \partial(|g_{1,1} * f|_{\mathrm{sqf}}) \in \mathfrak{m}\mathbb{F}.$$
In the last step we used that $\mathbb{F}$ is minimal, so $\partial \mathbb{F} \subseteq \mathfrak{m}\mathbb{F}$. Now, $\partial f$ has a smaller homological degree than $f$, so it can be decomposed into a product by induction.

Moreover, note that $r \in \mathfrak{m}\mathbb{F}$ implies that $r$ is an $S$-linear combination of elements of $\mathbb{F}$ of strictly smaller multidegree than $f$. So these elements can be decomposed by induction as well. To conclude that $f$ can be decomposed as claimed, we note the following fact: Whenever for $f' \in \mathbb{F}$ and a monomial $m \in S$ it holds that $m|f'|_{\mathrm{sqf}}$ is squarefree, then it follows that $m|f'|_{\mathrm{sqf}} = |m|f'|_{\mathrm{sqf}}|_{\mathrm{sqf}} = |mf'|_{\mathrm{sqf}}$. □

The preceding result does not hold if $I$ not squarefree.

**Example 3.5.** Consider the multiplication constructed in Example 3.3. Under that multiplication, $\mathbf{g}_{bc}$ cannot be written as a product of elements of degree one.

We will see in Section 6 below that every monomial ideal admits at least one multiplication which satisfies the conclusion of Proposition 3.4. This is no longer true for more general ideals. Indeed, if $I \subseteq S$ is a homogeneous ideal with three generators satisfying the conclusion of Proposition 3.4, then the third Betti number is bounded by the number of non-associative, graded-commutative products of three elements. On the other hand, it is a celebrated result of Bruns [Bru76] that any free resolution can be modified to yield a resolution of an ideal with three generators, so there does not exist a global bound on the third Betti number.

The preceding proposition leads to the following structure theorem for DGA structures on $\mathbb{F}$:

**Theorem 3.6.** *Let $I \subseteq S$ be a squarefree monomial ideal. Assume that the minimal resolution $\mathbb{F}$ of $S/I$ is a multigraded DGA.*

*Then there exists a DG-ideal $J \subseteq \mathbb{T}$ in the Taylor resolution of $S/I$, such that $\mathbb{F} \cong \mathbb{T}/J$ as multigraded DGAs over $S$.*

*Proof.* Let $Q := S[x_1^{-1}, \ldots, x_n^{-1}]$ denote the Laurent polynomial ring. As $\mathbb{T} \otimes_S Q$ is an exterior algebra over $Q$, we can define a map $\varphi : \mathbb{T} \otimes_S Q \to \mathbb{F} \otimes_S Q$ of DGAs by setting $\varphi(\mathbf{g}_m) := \mathbf{g}_m$, where the latter is interpreted as element of the algebraic Scarf complex of $S/I$.

We need to show that $\varphi$ restricts to a map $\varphi : \mathbb{T} \to \mathbb{F}$, where we consider $\mathbb{T}$ and $\mathbb{F}$ as subalgebras of $\mathbb{T} \otimes Q$ and $\mathbb{F} \otimes Q$ in the natural way. For this, consider a subset $A = \{m_1 < m_2 < \ldots < m_s\} \subseteq G(I)$. It holds that

$$\mathbf{g}_A = \frac{\operatorname{lcm}(m_1, \ldots, m_s)}{m_1 \cdots m_s} \mathbf{g}_{m_1} * \mathbf{g}_{m_2} * \cdots * \mathbf{g}_{m_s},$$

and, consequently, that

$$\varphi(\mathbf{g}_A) = \frac{\operatorname{lcm}(m_1, \ldots, m_s)}{m_1 \cdots m_s} \varphi(\mathbf{g}_{m_1}) * \varphi(\mathbf{g}_{m_2}) * \cdots * \varphi(\mathbf{g}_{m_s})$$
$$= |\varphi(\mathbf{g}_{m_1}) * \varphi(\mathbf{g}_{m_2}) * \cdots * \varphi(\mathbf{g}_{m_s})|_{\text{sqf}}.$$

The last expression is clearly contained in $\mathbb{F}$, and as $\mathbb{T}$ is generated by the elements $(\mathbf{g}_A)_{A \subseteq G(I)}$ it follows that $\varphi$ maps $\mathbb{T}$ to $\mathbb{F}$.

It remains to show that the restriction of $\varphi$ to $\mathbb{T}$ is surjective onto $\mathbb{F}$. But this is clear from Proposition 3.4, because $\varphi$ is surjective in homological degree 1 and it commutes with taking the squarefree part (the latter holds for any morphism which is homogeneous with respect to the multigrading). $\square$

The following example shows that one cannot omit the assumption of squarefreeness in the preceding Theorem.

**Example 3.7.** We continue with Example 3.3 and consider the DGAs structure on $\mathbb{F}$ constructed in that example. Let $\varphi : \mathbb{T} \otimes Q \to \mathbb{F} \otimes Q$ the map from the proof of Theorem 3.6. Then

$$\varphi(\mathbf{g}_{bc}) = \varphi(\frac{1}{x}\mathbf{g}_b * \mathbf{g}_c) = \frac{1}{x}\varphi(\mathbf{g}_b) * \varphi(\mathbf{g}_c) = -\frac{z}{x}\mathbf{g}_{ab} + \frac{y}{x}\mathbf{g}_{ac},$$

which is not contained in $\mathbb{F}$. So $\varphi$ does not restrict to a map $\mathbb{T} \to \mathbb{F}$ in this case. In fact, $\mathbb{F}$ is not the image of $\mathbb{T}$ under any map of DGAs, because $\mathbb{F} \otimes Q$ is not generated in homological degree 1.

The next example illustrates how one can use Theorem 3.6 to prove that a given ideal does not admit a (multigraded) minimal DGA resolution.

**Example 3.8.** Let $I \subseteq S = \mathbb{Q}[x_1, x_2, x_3, x_4, x_5, x_6]$ be the ideal generated by $a := x_1 x_2, b := x_2 x_3, c := x_3 x_4, d := x_4 x_5$ and $e := x_5 x_6$. This is a well-known example by Avramov which does not admit a minimal DGA resolution [Avr81, Example I].

Assume for the contrary that its minimal free resolution $\mathbb{F}$ admits a DGA structure. Then by Theorem 3.6 there exists a DG-ideal $J \subseteq \mathbb{T}$ in its Taylor resolution such that $\mathbb{F} = \mathbb{T}/J$. We have that $\beta^S_{3,(1,1,1,1,0,0)}(S/I) = 0$ and hence $\mathbf{g}_{abc} \in J$. As $J$ is a DG-ideal, it also holds that $(\partial \mathbf{g}_{abc}) * \mathbf{g}_e \in J$. Moreover, $\mathbf{g}_{abc} * \mathbf{g}_d = x_3 \mathbf{g}_{abcd} \in J$. As $\mathbb{T}/J$ is free, $J$ is saturated with respect to the variables and hence also $\mathbf{g}_{abcd} \in J$. Similarly, $\beta^S_{3,(0,0,1,1,1,1)}(S/I) = 0$ implies that $\mathbf{g}_a * (\partial \mathbf{g}_{cde}) \in J$ and $\mathbf{g}_{bcde} \in J$. It follows that

$$f := (\partial \mathbf{g}_{abc}) * \mathbf{g}_e - \mathbf{g}_a * (\partial \mathbf{g}_{cde}) - x_1 \partial \mathbf{g}_{bcde} - x_4 \partial \mathbf{g}_{abcd}$$
$$= x_2 x_3^2 \mathbf{g}_{abe} - x_2^2 x_3 \mathbf{g}_{ade} + x_1 x_2 x_3 \mathbf{g}_{bde} - x_2 x_3 x_4 \mathbf{g}_{abd}$$

is contained in $J$. On the other hand, the sets $\{a,b,e\}, \{a,d,e\}, \{b,d,e\}$ and $\{a,b,d\}$ lie in the Scarf complex of $I$, so the images of the corresponding generators of $\mathbb{T}$ are part of a basis of $\mathbb{F}$. But this is a contradiction, because $f$ is a relation between these elements.

## 4. Consequences of the existence

4.1. **Betti numbers.** In this section, we prove the following characterization of the possible Betti vectors of squarefree monomial ideals admitting minimal DGA resolutions.

**Theorem 4.1.** *Let $f = (1, f_1, f_2, \dots) \in \mathbb{N}^\nu$ be a finite sequence of natural numbers. Then the following conditions are equivalent:*
   *(1) There exists a squarefree monomial ideal $I$ in some polynomial ring $S$ whose minimal free resolution is a DGA, such that $\beta^S_i(S/I) = f_i$ for all $i$.*
   *(2) $f$ is the $f$-vector of a simplicial complex $\Delta$ which is a cone.*

Recall that the $f$-vectors of simplicial complexes are characterized by the Kruskal-Katona theorem [Sta96, Theorem 2.1]. From this one can derive an explicit characterization of the $f$-vectors of cones, see also [Kal85]. For the implication "1) $\Rightarrow$ 2)" we use the following general result.

**Proposition 4.2.** *Let $A = \bigoplus_{k \geq 0} A_k$ be a DGA over a field $\Bbbk$. Assume that $A_0 = \Bbbk$, that $A$ is generated in degree $1$, that $\dim_\Bbbk A_1 \leq \infty$ and finally that the differential of $A$ is not identically zero. Then the Hilbert function of $A$ equals the $f$-vector of a simplicial complex which is a cone.*

Aramova, Herzog and Hibi [AHH97, Theorem 4.1] characterize the Hilbert functions of graded-commutative algebras which are generated in degree one. These are indeed the same as $f$-vectors of simplicial complexes. Our result extends this, as we show showing that the existence of a non-zero differential yields a further restriction to the Hilbert series.

*Proof.* Let $C \subseteq A$ be the subalgebra of cycles. We claim that $C$ is also generated in degree 1 and that $A \cong C \oplus C[-1]$ as $\Bbbk$-vector space, where $C[-1]$ denotes the vector space $C$ with degrees shifted by one. These two claims are sufficient to prove our result. Indeed, by the above mentioned [AHH97, Theorem 4.1], the Hilbert function of $C$ is the $f$-vector of some simplicial complex $\Delta$. An elementary computation shows that in this case the Hilbert function of $A \cong C \oplus C[-1]$ is the $f$-vector of the cone over $\Delta$.

Now we turn to the proof of our claims. Let $C_1 \subseteq A_1$ be the space of cycles in $A_1$ and write $C'$ for the algebra generated by them. By assumption, the differential $\partial$ of $A$ does not vanish identically. As $A$ is generated in degree one, the Leibniz rule implies that the first component $\partial_1 : A_1 \to A_0 = \Bbbk$ is nonzero, hence there exists an element $\hat{a} \in A_1$ with $\partial \hat{a} = 1$ and (thus) $A_1 = C_1 \oplus \Bbbk \hat{a}$. The latter implies that $A \cong C' + \hat{a}C'$.

Hence any element $a \in A$ can be written as $a = c_1 + \hat{a}c_2$ with $c_1, c_2 \in C'$. In particular, if $a \in C$, then $0 = \partial a = c_2$, so $a = c_1 \in C'$ and thus $C = C'$ is generated in degree one.

For the second claim, we first note that $A \cong C \oplus \hat{a}C$. Indeed, for any $\hat{a}c \in C \cap \hat{a}C$ it holds that $0 = \partial(\hat{a}c) = c$. Finally, it remains to show that $\hat{a}C \cong C(-1)$. For this it is sufficient to show that the multiplication by $\hat{a}$ is injective on $C$. But this is clear, because $\partial(\hat{a}c) = c$ for any $c \in C$, so the differential provides a left-inverse. $\square$

*Proof of Theorem 4.1.* 1) $\Rightarrow$ 2) Let $Q$ be the field of fractions of $S$ and let $\mathbb{F}$ the minimal free resolution of $I$. It follows from Proposition 3.4 that $\mathbb{F} \otimes_S Q$ is generated in degree 1 as $Q$-algebra. Hence the claim is immediate from Proposition 4.2, applied to $\mathbb{F} \otimes_S Q$.

2) $\Rightarrow$ 1) Set $k := f_1$. Let $\Delta \subseteq 2^k$ be a simplicial complex which is a cone with apex 1 and whose $f$-vector equals $f$. We can find a monomial ideal $I$ in some polynomial ring $S$, such that the lcm-lattice of $I$ equals the face lattice of $\Delta$, augmented by a maximal element (cf. [IKMF17, Theorem 3.4], [Map13]). As this property is invariant under polarization, we may as well assume that $I$ is squarefree.

By construction, $\Delta$ is the Scarf complex of $I$. As $\Delta$ is acyclic, it follows that the algebraic Scarf complex is acyclic as well. But the latter is a subcomplex of the minimal free resolution of $S/I$, and thus the two coincide. Hence the $f$-vector of $\Delta$ equals the Betti vector of $I$.

It remains to show that the minimal free resolution of $I$ admits a DGA structure. Let $\mathbb{T}$ be Taylor resolution of $I$. We identify the set of generators of $I$ with the set $[k]$ of vertices of $\Delta$. To construct the minimal free resolution of $S/I$ we use discrete Morse theory [JW09; Skö06]. Let us recall the relevant aspects of this technique. Let $G$ be the directed graph with vertex set $2^{[k]}$ and edges from $V$ to $W$ whenever the coefficient of $\mathbf{g}_V$ in the expansion of $\partial \mathbf{g}_W$ is nonzero. A *Morse matching* $\mathcal{M}$ is a collection of edges $(V, W)$ in this graph with the following properties:

(a) $\mathcal{M}$ is a matching, i.e. the edges are disjoint,
(b) reversing all edges of $\mathcal{M}$ in $G$ results in an acyclic graph, and
(c) $\deg \mathbf{g}_V = \deg \mathbf{g}_W$ for all $(V, W) \in \mathcal{M}$.

Given such a Morse matching, define

$$J_\mathcal{M} := \mathrm{span}_\Bbbk \{\mathbf{g}_W, \partial \mathbf{g}_W : \exists V \in 2^{[k]}, (V, W) \in \mathcal{M}\}.$$

Then the complex $\mathbb{T}^\mathcal{M} := \mathbb{T}/J_\mathcal{M}$ is a free resolution of $S/I$, and it has a $\Bbbk$-basis is given by the images of $\mathbf{g}_W$ for those $W \in 2^{[k]}$ which do not appear in the matching.

In our situation we define the Morse matching:

$$\mathcal{M} := \{(\mathbf{g}_{W \setminus \{1\}}, \mathbf{g}_W) : W \subset [k], W \notin \Delta, 1 \in W\}$$

Recall that we assumed 1 to be the apex of $\Delta$, so $W \notin \Delta$ and $1 \in W$ imply that $W \setminus \{1\} \notin \Delta$. This is clearly a matching, and it is easy to see that it satisfies property (b). For property (c), recall that the lcm lattice $L_I$ of $I$ is $\Delta$, together with an additional maximal element. We are only matching generators which are not in $\Delta$, so the degrees of all those generators correspond to the maximal element of $L_I$, and thus they are equal.

Finally, the unmatched generators correspond exactly to the faces of $\Delta$. As we know that the Betti vector of $I$ equals the $f$-vector of $\Delta$, we can conclude that $\mathbb{T}^\mathcal{M} := \mathbb{T}/J_\mathcal{M}$ is a *minimal* free resolution of $S/I$. On the other hand, by our choice of the matching $\mathcal{M}$, $J_\mathcal{M}$ is in fact a DG-ideal and hence $\mathbb{T}^\mathcal{M}$ is a DG-algebra. □

We give an example to show that the assumptions of $I$ being squarefree and admitting a minimal DGA resolution are both necessary in Theorem 4.1.

**Example 4.3.** Consider the ideal $I := \langle x_1x_2, x_2x_3, x_3x_4, x_4x_5, x_5x_6, x_6x_1 \rangle \subseteq \mathbb{Q}[x_1, \ldots, x_6]$. Using `Macaulay2` [GS], one can compute that its total Betti numbers are $(1, 6, 9, 6, 2)$. We claim that this is not the $f$-vector of any simplicial complex. Indeed, such a simplicial complex would have two 3-simplices. As each 3-simplex has four 2-faces, and the two 3-simplices can only share one 2-face, we would need a total of at least seven 2-faces, but we have only six.

On the other hand, by Theorem 2.1 there exists a monomial $m \in S$ such that $mI$ admits a minimal DGA resolution. However, $mI$ is not squarefree, and indeed its Betti vector still does not satisfy the conclusion of Theorem 4.1, as it coincides with the Betti vector of $I$.

4.2. **Subadditivity of syzygies.** For a monomial ideal $I \subseteq S$ and $0 \leq i \leq \operatorname{pdim} S/I$ we define
$$t_i := \max\{j : \beta^S_{i,j}(S/I) \neq 0\}.$$
We say that the syzygies of $I$ are *subadditive* when $t_b \leq t_a + t_{b-a}$ for all $a, b$ such that $1 \leq a < b \leq \operatorname{pdim} S/I$. Not every ideal has this property [ACI15], but among monomial ideals no counterexample is known. Here, we show that for squarefree monomial ideals, the existence of a DGA structure on the minimal free resolution implies the subadditivity of syzygies.

**Proposition 4.4.** *If $I \subseteq S$ is a squarefree monomial ideal which admits a minimal DGA resolution $\mathbb{F}$, then its syzygies are subadditive.*

The following lemma is needed for the proof of this result.

**Lemma 4.5.** *Let $f, g \in \mathbb{F}$ be homogeneous. If $|f * g|_{\operatorname{sqf}} \notin \mathfrak{m}\mathbb{F}$, then there exist $f', g' \in \mathbb{F} \setminus \mathfrak{m}\mathbb{F}$ of the same homological degrees as $f$ and $g$, respectively, such that $\deg(f) \vee \deg(g) = \deg(f') \vee \deg(g')$.*

*Proof.* By choosing a basis of $\mathbb{F}$, we may write $f = \sum_i \lambda_i c_i f_i$ and $g = \sum_j \mu_j d_j g_j$, with $\lambda_i, \mu_i \in \Bbbk$, $c_i, d_j \in S$ monomials and $f_i, g_j \in \mathbb{F} \setminus \mathfrak{m}\mathbb{F}$. By our assumption, $|f * g|_{\operatorname{sqf}} = \sum_{i,j} \lambda_i \mu_j |c_i d_j f_i * g_j|_{\operatorname{sqf}} \notin \mathfrak{m}\mathbb{F}$. So at least one summand is not contained in $\mathfrak{m}\mathbb{F}$, say $|c_1 d_1 f_1 * g_1|_{\operatorname{sqf}}$. It follows easily from the definition of the squarefree part that there exists a monomial $m \in S$ such that $|c_1 d_1 f_1 * g_1|_{\operatorname{sqf}} = m|f_1 * g_1|_{\operatorname{sqf}}$. But this does not lie in $\mathfrak{m}\mathbb{F}$, so we can conclude that $m = 1$ and the claim follows with $f' := f_1$ and $g' := g_1$. □

*Proof of Proposition 4.4.* Let $f \in \mathbb{F}_b$ be an element of total degree $t_b$ and we may assume that $f \notin \mathfrak{m}\mathbb{F}$. By Proposition 3.4, we can write $f$ as
$$f = \sum_j |g_{1,j} * (g_{2,j} * (\cdots (g_{i-1,j} * g_{i,j})\cdots))|_{\text{sqf}},$$
for $g_{\ell,j} \in \mathbb{F}_1$ for $1 \leq \ell \leq i$ and all $j$. As $f \notin \mathfrak{m}\mathbb{F}$, the same holds for at least one summand, say for $j = 1$.

Using that the multiplication on $\mathbb{F}$ is associative, we may rewrite that summand as $|(g_{1,1}\cdots g_{a,1}) * (g_{a+1,1}\cdots g_{b,1})|_{\text{sqf}}$. Now we apply Lemma 4.5 to this product to obtain $f' \in \mathbb{F}_a$ and $g' \in \mathbb{F}_{b-a}$ which are both not contained in $\mathfrak{m}\mathbb{F}$ and satisfy $\deg f = \deg(f') \vee \deg(g')$. Thus we conclude that
$$t_b = \deg(f) = \deg(f') \vee \deg(g') \leq \deg(f') + \deg(g') \leq t_a + t_{b-a}$$
$\square$

If we do not assume that $\mathbb{F}$ admits a DGA structure, then our methods still suffice to show the case $a = 1$. This was first shown by Herzog and Srinivasan in [HS15, Corollary 4] using a different method.

**Corollary 4.6.** *For a monomial ideal $I \subseteq S$, it holds that $t_i \leq t_1 + t_{i-1}$ for all $2 \leq i \leq \text{pdim}\, S/I$.*

*Proof.* We may replace $I$ by its polarization and so assume it is squarefree. Further, choose any multiplication on $\mathbb{F}$. Now we just note that the proof of Proposition 4.4 does not require the multiplication to be associative if $a = 1$. $\square$

## 5. Strongly generic ideals

A monomial ideal is called *strongly generic* if no variable appears with the same nonzero exponent in two distinct generators of $I$. This class of ideals was introduced and studied by Bayer, Peeva and Sturmfels [BPS98] [1], see also Chapter 6 of [MS05]. Strongly generic ideals have a number of desirable properties, and in particular their minimal free resolution is given by the algebraic Scarf complex [BPS98, Theorem 3.2]. However, this condition is not sufficient to ensure that the ideal admits a minimal DGA resolution.

**Theorem 5.1.** *The ideal $\langle x^2, xy, y^2z^2, zw, w^2 \rangle \subseteq \Bbbk[x, y, z, w]$ is strongly generic, but its minimal free resolution does not admit the structure of a DGA, whose multiplication respects the standard $\mathbb{Z}$-grading.*

This example is a variation of Example 3.8. Note that the theorem not only states that there is no multiplication respecting the multigrading, but that there is even no multiplication satisfying the weaker assumption of respecting the standard $\mathbb{Z}$-grading. The proof of Theorem 5.1 is very technical, so we postpone it and discuss the result first.

**Remark 5.2.** Corollary 3.6 in [BPS98] states that strongly generic monomial ideals admit a minimal DGA resolution. I was informed by I. Peeva that the authors of that article are aware that the claim does not hold as stated and that the statement and the proof were supposed to include an additional combinatorial condition.

---
[1] In [BPS98] these ideals are called *generic*. In the later article [MSY00] a weaker notion of genericity was introduced and the original one is called strongly generic since then.

Theorem 5.1 implies that various kinds of resolutions do not admit a DGA structure in general.

**Corollary 5.3.**   *(1) There exists a monomial ideal whose hull resolution [BS98] does not admit a DGA structure.*
  *(2) There exists a monomial ideal whose Lyubeznik resolution [Lyu88] does not admit a DGA structure.*
  *(3) There exists a monomial ideal whose minimal free resolution is supported on a simplicial complex (and in particular cellular) [BS98], but it does not admit a DGA structure.*

*Proof.* The ideal of Theorem 5.1 provides the counterexample in all three cases. This ideal is generic, so its minimal free resolution is supported in its Scarf complex and it coincides with its hull resolution [MS05, Theorem 6.13]. The Lyubeznik resolution depends on the choice of a total order on the generators. One can check that the Lyubeznik resolution with respect to the order $xy \prec zw \prec x^2 \prec y^2z^2 \prec w^2$ is minimal. □

In fact, the minimal free resolution of Avramov's Example 3.8 can also be given the structure of a cellular resolution, though not with a simplicial complex as support.

*Proof of Theorem 5.1.* In this proof we use both the multigrading and the standard $\mathbb{Z}$-grading, so we refer to the latter as the total degree.

Let $S = \Bbbk[x, y, z, w]$ and let $I$ be the ideal of the statement. We give names to its generators as follows:

$$a := x^2, \ b := xy, \ c := y^2z^2, \ d := zw, \ e := w^2$$

The Scarf complex $\Delta_I$ is the simplicial complex with the facets $\{b, c, d\}$ and $\{a, b, d, e\}$. As mentioned above, the minimal free resolution $\mathbb{F}$ of $S/I$ is given by its algebraic Scarf complex. The generators of $\mathbb{F}$ are listed in Table 2, together with their total degrees.

Assume that $\mathbb{F}$ admits an associative multiplication which respects the $\mathbb{Z}$-grading, which we denote by $*$. The map $* : \mathbb{F} \otimes \mathbb{F} \to \mathbb{F}$ can be decomposed into its homogeneous components with respect to the $\mathbb{Z}^4$-grading. We write $\times$ for the component of multidegree $(0, 0, 0, 0)$. Because the differential on $\mathbb{F}$ respects the multigrading, $\times$ also satisfies the Leibniz rule and lifts the multiplication on $S/I$. Hence, by the uniqueness of the multiplication, there exists a map $\sigma : \mathbb{F} \otimes \mathbb{F} \to \mathbb{F}[1]$ such that

$$a * b = a \times b + \delta\sigma(a, b)$$

for $a, b \in \mathbb{F}$, where $\delta\sigma = \partial \circ \sigma + \sigma \circ \partial$. As we assumed that $*$ respects the total degree, it follows that $\sigma$ is homogeneous with respect to it. However, in general $\sigma$ is not homogeneous with respect to the multigrading. In fact, the component of $\sigma$ in multidegree $(0, 0, 0, 0)$ vanishes, because $*$ and $\times$ coincide in this degree. We proceed by proving a series of claims.

**Claim 1.** *For $i, j, k \in \{a, b, d, e\}$, all products of the form $\mathbf{g}_i \times \mathbf{g}_j$ and $\mathbf{g}_i \times \mathbf{g}_{jk}$ coincide with the corresponding products from the Taylor resolution.*

| | |
|---:|:---:|
| $\mathbf{g}_a, \mathbf{g}_b, \mathbf{g}_d, \mathbf{g}_e$ | 2 |
| $\mathbf{g}_c$ | 4 |
| $\mathbf{g}_{ab}, \mathbf{g}_{de}$ | 3 |
| $\mathbf{g}_{ad}, \mathbf{g}_{ae}, \mathbf{g}_{bd}, \mathbf{g}_{be}$ | 4 |
| $\mathbf{g}_{bc}, \mathbf{g}_{cd}$ | 5 |
| $\mathbf{g}_{ijk} \neq \mathbf{g}_{bcd}$ | 5 |
| $\mathbf{g}_{bcd}$ | 6 |
| $\mathbf{g}_{abde}$ | 6 |

TABLE 2. The degrees of the generators of $\mathbb{F}$.

If $I$ were squarefree, then this would follow from Proposition 3.1. However, as $I$ is not squarefree, we need to verify this. We start in homological degree 2. For any $i, j \in \{a, b, d, e\}, i \neq j$, one can verify that the product $\mathbf{g}_i \times \mathbf{g}_j$ is a multiple of $\gcd(i,j)\mathbf{g}_{ij}$, because there is no other element of the correct multidegree. Then the Leibniz rule implies that $\mathbf{g}_i \times \mathbf{g}_j = \pm \gcd(i,j)\mathbf{g}_{ij}$. Moreover $\mathbf{g}_i \times \mathbf{g}_i = 0$, because $*$ and thus $\times$ was assumed to be graded-commutative.

Next, $\mathbf{g}_i \times \mathbf{g}_{ij} = 0$ for all $i, j$ because it is clearly a cycle, and the only boundary of homological degree 3 is $\partial \mathbf{g}_{abde}$, which has a different multidegree. Moreover, for pairwise distinct $i, j, k \in \{a, b, d, e\}$, it follows again by inspection that $\mathbf{g}_i \times \mathbf{g}_{jk}$ is a multiple of $\gcd(i, \text{lcm}(j,k))\mathbf{g}_{ijk}$ because this is the only element with the required multidegree, and again the Leibniz rule determines the coefficient to be $\pm 1$.

In fact, the other products in the subalgebra generated by $\mathbf{g}_a, \mathbf{g}_b, \mathbf{g}_d$ and $\mathbf{g}_d$ also coincide with the product from the Taylor resolution, but we do not need this.

**Claim 2.** *It holds that* $(\mathbf{g}_a \times \mathbf{g}_c) \times \mathbf{g}_e - \mathbf{g}_a \times (\mathbf{g}_c \times \mathbf{g}_e) = yz\partial \mathbf{g}_{abde}$

By considering the multigrading and using the Leibniz rule as above, one can show that $\mathbf{g}_a \times \mathbf{g}_c = yz^2\mathbf{g}_{ab} + x\mathbf{g}_{bc}$ and $\mathbf{g}_c \times \mathbf{g}_e = w\mathbf{g}_{cd} + y^2z\mathbf{g}_{de}$. Using this and Claim 1, one can compute that

$$\partial(\mathbf{g}_{bc} \times \mathbf{g}_e) = \partial(w\mathbf{g}_{bcd} + yz\mathbf{g}_{bde}).$$

This implies that $\mathbf{g}_{bc} \times \mathbf{g}_e = w\mathbf{g}_{bcd} + yz\mathbf{g}_{bde} + \partial(f)$ for some $f \in \mathbb{F}_4$. But the $x$-component of the multidegree of every element in $\mathbb{F}_4$ is at least 2, hence $f = 0$. Now we can compute that

$$(\mathbf{g}_a \times \mathbf{g}_c) \times \mathbf{g}_e = yz^2\mathbf{g}_{ab} \times \mathbf{g}_e + x\mathbf{g}_{bc} \times \mathbf{g}_e = yz^2\mathbf{g}_{abe} + xw\mathbf{g}_{bcd} + xyz\mathbf{g}_{bde}.$$

Using essentially the same reasoning one can also show that

$$\mathbf{g}_a \times (\mathbf{g}_c \times \mathbf{g}_e) = y^2z\mathbf{g}_{ade} + xw\mathbf{g}_{bcd} + yzw\mathbf{g}_{abd}.$$

Hence it follows that $(\mathbf{g}_a \times \mathbf{g}_c) \times \mathbf{g}_e - \mathbf{g}_a \times (\mathbf{g}_c \times \mathbf{g}_e) = yz\partial \mathbf{g}_{abde}$.

**Claim 3.** *For any $f \in \mathbb{F}$ it holds that $\sigma(1_\Bbbk, f) = \sigma(f, 1_\Bbbk) = 0$, where $1_\Bbbk$ denotes the "1" in $\Bbbk$.*

For $f \in \mathbb{F}$ it holds that

$$f = 1_\Bbbk * f = 1_\Bbbk \times f + \partial\sigma(1_\Bbbk, f) + \underbrace{\sigma(\partial 1_\Bbbk, f)}_{=0} + \sigma(1_\Bbbk, \partial f)$$

and hence $\partial\sigma(1_\Bbbk, f) + \sigma(1_\Bbbk, \partial f) = 0$.

We proceed now by induction on $|f|$. Consider the case $|f| = 0$, i.e., $f = 1_\Bbbk$. Then $\sigma(1_\Bbbk, 1_\Bbbk)$ has multidegree $(0, 0, 0, 0)$ and homological degree 1, which is only possible if $\sigma(1_\Bbbk, 1_\Bbbk) = 0$.

Assume now that $|f| = i$ and $\sigma(1_\Bbbk, g) = 0$ for all $g \in \mathbb{F}$ with $|g| < i$. Then $\partial \sigma(1_\Bbbk, f) = -\sigma(1_\Bbbk, \partial f) = 0$ by the induction hypothesis. As $\mathbb{F}$ is acyclic, there exist an $h \in \mathbb{F}$ with $|h| = |f| + 2$ such that $\partial h = \sigma(1_\Bbbk, f)$. We may assume that $f$ is a generator of $\mathbb{F}$. Now, inspection of Table 2 shows that for no generator $f$ of $\mathbb{F}$ there exist an element $h$ with $\deg h = \deg \sigma(1_\Bbbk, f) = \deg f$ and $|h| = |f| + 2$. Hence $\sigma(1_\Bbbk, f) = 0$.

**Claim 4.** *It holds that $\sigma(\mathbf{g}_i, \mathbf{g}_j) = 0$ for $i, j \neq c$.*

If $\sigma(\mathbf{g}_i, \mathbf{g}_j) \neq 0$ for $i, j \neq c$, then it is an element of total degree 4 and homological degree 3. But by Table 2, no such element exists.

**Claim 5.** *It holds that $\mathbf{g}_i * \mathbf{g}_j = \mathbf{g}_i \times \mathbf{g}_j$ for $i, j \neq c$.*

In general we have that $\mathbf{g}_i * \mathbf{g}_j = \mathbf{g}_i \times \mathbf{g}_j + \partial \sigma(\mathbf{g}_i, \mathbf{g}_j) + \sigma(\partial \mathbf{g}_i, \mathbf{g}_j) + \sigma(\mathbf{g}_i, \partial \mathbf{g}_j)$. All terms involving $\sigma$ vanish because of Claim 3 and Claim 4.

**Claim 6.** *For $i, j, k \in \{a, b, d, e\}$, it holds that $\sigma(\mathbf{g}_i, \mathbf{g}_{jk}) = 0$.*

First, assume that $i = j$ or $i = k$. We only consider the case $i = j$, the other case is similar. For this, we compute that

$$0 = (\mathbf{g}_i * \mathbf{g}_i) * \mathbf{g}_j = \mathbf{g}_i * (\mathbf{g}_i * \mathbf{g}_j) \stackrel{\text{Cl. } 5}{=} \mathbf{g}_i * (\mathbf{g}_i \times \mathbf{g}_j) \stackrel{\text{Cl. } 1}{=} \pm \mu \mathbf{g}_i * \mathbf{g}_{ij}$$
$$= \pm \mu (\underbrace{\mathbf{g}_i \times \mathbf{g}_{ij}}_{=0} + \delta \sigma(\mathbf{g}_i, \mathbf{g}_{ij})) \stackrel{\text{Cl. } 3}{=} \pm \mu \partial \sigma(\mathbf{g}_i, \mathbf{g}_{ij}),$$

where $\mu = \gcd(i, j)$. The differential is injective in homological degree 4, hence $\partial \sigma(\mathbf{g}_i, \mathbf{g}_{ij}) = 0$ implies that $\sigma(\mathbf{g}_i, \mathbf{g}_{ij}) = 0$.

Now we consider the case that $i, j, k$ are all different. Then $\{i, j, k\}$ contains either $a$ and $b$, or $d$ and $e$. We only consider the first case, the other one is similar. If $\{j, k\} = \{a, b\}$ then $\sigma(\mathbf{g}_i, \mathbf{g}_{ab}) = 0$ holds for degree reasons, see Table 2. So assume that $i \in \{a, b\}$. Again, we assume without loss of generality that $i = a$. Then either $j$ or $k$ equals $b$, say $j$. We compute that

$$(\mathbf{g}_a * \mathbf{g}_b) * \mathbf{g}_k \stackrel{\text{Cl. } 5}{=} (\mathbf{g}_a \times \mathbf{g}_b) * \mathbf{g}_k = x \mathbf{g}_{ab} * \mathbf{g}_k = x \mathbf{g}_{ab} \times \mathbf{g}_k + x \delta \sigma(\mathbf{g}_{ab}, \mathbf{g}_k)$$
$$\stackrel{(*)}{=} x \mathbf{g}_{abk} + x \delta \sigma(\mathbf{g}_{ab}, \mathbf{g}_k)$$
$$\stackrel{\text{Cl. } 4}{=} x \mathbf{g}_{abk} + x \partial \underbrace{\sigma(\mathbf{g}_{ab}, \mathbf{g}_k)}_{=0},$$

where in $(*)$ we used Claim 1 and that $\gcd(\text{lcm}(a, b), k) = 1$ for $k \in \{d, e\}$. On the other hand, a similar computation shows that $\mathbf{g}_a * (\mathbf{g}_b * \mathbf{g}_k) = x \mathbf{g}_{abk} + \partial \sigma(\mathbf{g}_a, \mathbf{g}_{bk})$, where we use that $\gcd(b, k) = 1$ and $\gcd(a, \text{lcm}(b, k)) = x$. The assumption that $*$ is associative implies that $\partial \sigma(\mathbf{g}_a, \mathbf{g}_{bk}) = 0$ and thus $\sigma(\mathbf{g}_a, \mathbf{g}_{bk})$ because $\partial$ is injective in homological degree 4.

**Claim 7.** *It holds that $\mathbf{g}_a * (\mathbf{g}_c * \mathbf{g}_e) - (\mathbf{g}_a * \mathbf{g}_c) * \mathbf{g}_e \neq 0$.*

Let us compute the associator of $\mathbf{g}_a, \mathbf{g}_c$ and $\mathbf{g}_e$:

$$\mathbf{g}_a * (\mathbf{g}_c * \mathbf{g}_e) - (\mathbf{g}_a * \mathbf{g}_c) * \mathbf{g}_e =$$
$$\mathbf{g}_a \times (\mathbf{g}_c \times \mathbf{g}_e) - (\mathbf{g}_a \times \mathbf{g}_c) \times \mathbf{g}_e +$$
$$\delta\sigma(\mathbf{g}_a, \mathbf{g}_c \times \mathbf{g}_e) - \delta\sigma(\mathbf{g}_a \times \mathbf{g}_c, \mathbf{g}_e) +$$
$$\mathbf{g}_a \times \delta\sigma(\mathbf{g}_c, \mathbf{g}_e) - \delta\sigma(\mathbf{g}_a, \mathbf{g}_c) \times \mathbf{g}_e +$$
$$\delta\sigma(\mathbf{g}_a, \delta\sigma(\mathbf{g}_c, \mathbf{g}_e)) - \delta\sigma(\delta\sigma(\mathbf{g}_a, \mathbf{g}_c), \mathbf{g}_e)$$

We know that $\mathbf{g}_a \times (\mathbf{g}_c \times \mathbf{g}_e) - (\mathbf{g}_a \times \mathbf{g}_c) \times \mathbf{g}_e = yz\partial\mathbf{g}_{abde}$ and we want to show that this term does not cancel against any other term. For this we need to consider the component of $\mathbf{g}_a * (\mathbf{g}_c * \mathbf{g}_e) - (\mathbf{g}_a * \mathbf{g}_c) * \mathbf{g}_e$ in the multidegree of $\mathbf{g}_a \times (\mathbf{g}_c \times \mathbf{g}_e)$. As $\sigma$ has no component in multidegree $(0,0,0,0)$, the summands having only one $\sigma$ cannot contribute. Thus it remains to consider the last two terms, namely $\delta\sigma(\mathbf{g}_a, \delta\sigma(\mathbf{g}_c, \mathbf{g}_e))$ and $\delta\sigma(\delta\sigma(\mathbf{g}_a, \mathbf{g}_c), \mathbf{g}_e)$.

By Claim 3 it holds that $\delta\sigma(\mathbf{g}_a, \delta\sigma(\mathbf{g}_c, \mathbf{g}_e)) = \partial\sigma(\mathbf{g}_a, \partial\sigma(\mathbf{g}_c, \mathbf{g}_e))$. Further, $\sigma(\mathbf{g}_c, \mathbf{g}_e)$ has homological degree 3 and total degree 6, hence it is of the form $\sigma(\mathbf{g}_c, \mathbf{g}_e) = \lambda \mathbf{g}_{bcd} + r$ where $\lambda \in \Bbbk$ and $r$ is a combination of basis elements $\mathbf{g}_{ijk}$ with $i, j, k \neq c$. Now it holds that

$$\partial\sigma(\mathbf{g}_c, \mathbf{g}_e) = \lambda \partial \mathbf{g}_{bcd} + \partial r$$
$$= \lambda x \mathbf{g}_{cd} + \lambda w \mathbf{g}_{bc} + \text{terms not involving } c$$

Using Claim 6 it follows that

$$\partial\sigma(\mathbf{g}_a, \partial\sigma(\mathbf{g}_c, \mathbf{g}_e)) = \lambda x \partial\sigma(\mathbf{g}_a, \mathbf{g}_{cd}) + \lambda w \partial\sigma(\mathbf{g}_a, \mathbf{g}_{bc}) = \lambda(m_1 x + m_2 w)\partial\mathbf{g}_{abde},$$

for some $m_1, m_2 \in S$. Here we used that $\mathbf{g}_{abde}$ is the only generator in homological degree 4. These terms cannot cancel against $yz\partial\mathbf{g}_{abde}$, because they are multiples of $x$ or $w$, respectively. The same arguments apply to $\delta\sigma(\delta\sigma(\mathbf{g}_a, \mathbf{g}_c), \mathbf{g}_e)$. □

**Remark 5.4.** Unfortunately, it does not seem possible to use Avramov's obstructions from [Avr81] to prove the preceding theorem. At least, I could not find a suitable regular sequence $f_1, \ldots, f_r$ in the ideal to consider $R = S/(f_1, \ldots, f_s)$.

## 6. The lcm lattice and non-squarefree ideals

The *lcm lattice* $L_I$ of a monomial ideal is the lattice of all least common multiples of subsets of the minimal generators of $I$. The lcm lattice was introduced in [GPW99], where it is also proven that its isomorphism type determines the minimal free resolution of $I$. So one might be led to conjecture that it also determines the possible multiplication on $\mathbb{F}$. However, it is immediately clear from Theorem 2.1 that is not true. Nevertheless, we can identify a class of multiplications which is compatible with the lcm lattice:

**Definition 6.1.** A multiplication $*$ on $\mathbb{F}$ is called *supportive* if for all $a, b \in \mathbb{F}$ there exists a $c \in \mathbb{F}$ with $\deg c \leq (\deg a) \vee (\deg b)$, such that $a * b = mc$ for a monomial $m \in S$.

**Example 6.2.** The DGA structure on the Taylor resolution is supportive, but the multiplication constructed in Example 3.3 is not supportive.

For squarefree monomial ideals, being supportive is not a restriction.

**Lemma 6.3.** *If $I \subseteq S$ is a squarefree monomial ideal, then every multiplication on its minimal free resolution is supportive.*

*Proof.* Let $\mathbb{F}$ be the minimal free resolution of $S/I$ and let $\mathcal{B}$ be an $S$-basis for it. Every element of $\mathcal{B}$ has a squarefree degree, we may apply Lemma 1.4. It follows that for any $\mathbf{g}_1, \mathbf{g}_2 \in \mathcal{B}$ it holds that $\mathbf{g}_1 * \mathbf{g}_2 = mh$ for a monomial $m$ and an element $h \in \mathbb{F}$ with $\deg h \leq (\deg \mathbf{g}_1 + \deg \mathbf{g}_2) \wedge (1,\ldots,1) = \deg \mathbf{g}_1 \vee \deg \mathbf{g}_2$. $\square$

The next result shows that supportive multiplications are essentially determined by the isomorphism type of the lcm-lattice of $I$. This is our main motivation for introducing this notion.

**Proposition 6.4.** *Let $I$ and $I'$ be two monomial ideals in two polynomial rings $S$ and $S'$, respectively, whose lcm-lattices are isomorphic. Then every supportive multiplication on the minimal free resolution of $S/I$ can be relabeled to a supportive multiplication on the minimal free resolution of $S'/I'$.*

*Proof.* Let $\nu : L_I \to L_{I'}$ be the isomorphism. We recall the "relabeling" construction, which was introduced in [GPW99]. Fix a basis $\mathcal{B}$ of the minimal free resolution $\mathbb{F}$ of $S/I$. We express the differential $\partial$ of $\mathbb{F}$ in this basis:
$$\partial \mathbf{g} = \sum_{\mathbf{h} \in \mathcal{B}} c_{\mathbf{g}}^{\mathbf{h}} x^{\deg(\mathbf{g})-\deg(\mathbf{h})} \mathbf{h}$$
with $\mathbf{g} \in \mathcal{B}, c_{\mathbf{g}}^{\mathbf{h}} \in \Bbbk$. The *relabeled* resolution $\nu(\mathbb{F})$ is the free $S'$-module with the basis $\{\bar{\mathbf{g}} : \mathbf{g} \in \mathcal{B}\}$ and $\deg \bar{\mathbf{g}} = \nu(\deg \mathbf{g})$. Here, by $\bar{\mathbf{g}}$ we mean a new basis element. The differential on $\nu(\mathbb{F})$ is given by
$$\bar{\partial} \bar{\mathbf{g}} := \sum_{\mathbf{h} \in \mathcal{B}} c_{\mathbf{g}}^{\mathbf{h}} x^{\deg(\bar{\mathbf{g}})-\deg(\bar{\mathbf{h}})} \bar{\mathbf{h}}.$$
This is well-defined because $\deg \mathbf{h} \leq \deg \mathbf{g}$ implies that $\deg \bar{\mathbf{h}} = \nu(\deg \mathbf{h}) \leq \nu(\deg \mathbf{g}) = \deg \bar{\mathbf{g}}$. By [GPW99, Theorem 3.3], $\nu(\mathbb{F})$ is indeed a minimal free resolution of $S'/I'$. Our definition of the relabeling of the multiplication is analogous. If
$$\mathbf{g}_1 * \mathbf{g}_2 = \sum_{\mathbf{h} \in \mathcal{B}} c_{\mathbf{g}_1,\mathbf{g}_2}^{\mathbf{h}} x^{\deg(\mathbf{g}_1)+\deg(\mathbf{g}_2)-\deg(\mathbf{h})} \mathbf{h},$$
then we set
$$\bar{\mathbf{g}}_1 * \bar{\mathbf{g}}_2 := \sum_{\mathbf{h} \in \mathcal{B}} c_{\mathbf{g}_1,\mathbf{g}_2}^{\mathbf{h}} x^{\deg(\bar{\mathbf{g}}_1)+\deg(\bar{\mathbf{g}}_2)-\deg(\bar{\mathbf{h}})} \bar{\mathbf{h}}.$$
To see that this is well-defined, consider an $\mathbf{h} \in \mathcal{B}$ with $c_{\mathbf{g}_1,\mathbf{g}_2}^{\mathbf{h}} \neq 0$ and compute that

(1)
$$\deg \bar{\mathbf{g}}_1 + \deg \bar{\mathbf{g}}_2 \geq \deg \bar{\mathbf{g}}_1 \vee \deg \bar{\mathbf{g}}_2 = \nu(\deg \mathbf{g}_1) \vee \nu(\deg \mathbf{g}_2)$$
$$\stackrel{(*)}{=} \nu(\deg \mathbf{g}_1 \vee \deg \mathbf{g}_2) \stackrel{(**)}{\geq} \nu(\deg \mathbf{h}) = \deg \bar{\mathbf{h}}.$$

Here, the equality marked with $(*)$ holds because $\nu$ is an isomorphism of lattices. Moreover, the equality marked with $(**)$ holds because the multiplication is supportive.

It is clear that the relabeled multiplication is graded-commutative and satisfies the Leibniz rule, because the structure constants are the same as for the multiplication on $\mathbb{F}$. Finally, the computation (1) also shows that it is supportive. $\square$

The last two results have a number of immediate corollaries.

**Corollary 6.5.** *Every monomial ideal admits a supportive multiplication on its minimal free resolution.*

*Proof.* Choose any multiplication on the minimal free resolution of the polarization of $I$ and relabel it. □

**Corollary 6.6.** *Consider a supportive multiplication on $\mathbb{F}$. Then for every element $g \in \mathbb{F}$, there exists a monomial $m \in S$ such that $mg$ is a sum of products of elements of $\mathbb{F}_1$.*

*Proof.* If $I$ is squarefree, then this is immediate from Proposition 3.4. Moreover, the claim is clearly invariant under relabeling. Hence relabeling the supportive multiplication to the minimal free resolution of the polarization of $I$ yields the claim. □

**Corollary 6.7.** *Let $I$ and $I'$ be two monomial ideals with isomorphic lcm lattices and assume that $I$ is squarefree. If $I$ admits a minimal DGA resolution, then so does $I'$.*

*Proof.* Associativity is invariant under relabeling. □

We close the section by giving an example which illustrates that the last corollary cannot be strengthened by considering the Betti poset [CM14a; TV15; CM14b]. Recall that the *Betti poset* of a monomial ideal is the subposet of the lcm lattice consisting of those multidegrees in which are nonzero Betti numbers.

**Example 6.8.** The ideal $\langle x^2, xy, y^2z^2, zw, w^2 \rangle \subseteq \Bbbk[x,y,z,w]$ from Theorem 5.1 does not admit a minimal DGA resolution. As argued in the proof of that theorem, its minimal free resolution is given by its algebraic Scarf complex, so in particular its Betti poset is the face poset of its Scarf complex. However, that Scarf complex is a cone (with apex $b$ in the notation of that proof). Hence the construction in the proof of the second implication of Theorem 4.1 yields an example of a squarefree monomial ideal admitting a minimal DGA resolution but having the same Betti poset.

### Acknowledgments

The author is indebted to Volkmar Welker for the suggestion to look for a Macaulay-type theorem for DGA resolutions, which eventually led to Theorem 4.1. The author thanks the anonymous referees for several helpful comments.

Institut für Mathematik, Goethe-Universität, Frankfurt am Main, Germany
*E-mail address*: katthaen@math.uni-frankfurt.de